\documentclass[11pt, oneside]{article}
\usepackage{geometry}                		% See geometry.pdf to learn the layout options. There are lots.
\usepackage{graphicx}				% Use pdf, png, jpg, or eps\UTF{00C2}§ with pdflatex; use eps in DVI mode

\usepackage[usenames]{xcolor}		
\usepackage{amsmath,amsthm,amssymb}
\usepackage{amscd}

\usepackage{tikz}

\theoremstyle{definition}
\newtheorem{dfn}{Definition}

\theoremstyle{plain}
\newtheorem{theorem}[dfn]{Theorem}
\newtheorem{proposition}[dfn]{Proposition}

\newcommand{\Ag}{{\mathbb{A}_g}}
\newcommand{\bC}{\mathbb{C}}
\newcommand{\bR}{\mathbb{R}}
\newcommand{\bZ}{\mathbb{Z}}
\newcommand{\cH}{\mathcal{H}}
\newcommand{\Hom}{\mathrm{Hom}}
\newcommand{\Hg}{\mathfrak{H}_g}
\newcommand{\hP}{\widehat{\Phi}}
\newcommand{\Jac}{\mathrm{Jac}}
\newcommand{\Ker}{\mathrm{Ker}}
\newcommand{\Mg}{{\mathbb{M}_g}}
\newcommand{\la}{\lambda}
\newcommand{\om}{\widehat{\omega}}
\newcommand{\pa}{\partial}
\newcommand{\pab}{\overline{\partial}}
\newcommand{\smooth}{C^\infty}

\newcommand{\ib}{\overline{i}}
\newcommand{\jb}{\overline{j}}
\newcommand{\kb}{\overline{k}}
\newcommand{\lb}{\overline{l}}
\newcommand{\zb}{\overline{z}}

\title{A twisted invariant of a compact Riemann surface
\thanks
{%AMS subject classifications: Primary ??N05; Secondary ??F38, ??R15.
Keywords: moduli of compact Riemann surfaces, the first Mumford-Morita-Miller class, Kawazumi-Zhang invariant
}}

\date{}
\author{Nariya Kawazumi}

\begin{document}

\maketitle

\begin{abstract} We introduce a twisted version of the Kawazumi-Zhang invariant $a_g(C) = \varphi(C)$ of a compact Riemann surface $C$ of genus $g \geq 1$, and discuss how it is related to the first Mumford-Morita-Milller class $e_1 = \kappa_1$ on the moduli space of compact Riemann surfaces and the original Kawazumi-Zhang invariant.
\end{abstract}

\section{Introduction}
\label{sec:intro}
A compact Riemann surface is still a fascinating research subject even in the twenty-first century. In fact, various mathematicians and even physicists are now studying the Teichm\"uller space $\mathcal{T}_g$ of genus $g$, a universal covering space of the moduli space $\Mg$ of compact Riemann surfaces of genus $g$,
from their own various viewpoints. In this paper we study functions on the moduli space $\Mg$, i.e., invariant functions on the Teichm\"uller space $\mathcal{T}_g$ 
under the mapping class group action. They are regarded also as analytic invariants of a compact Riemann surface. There are no non-constant holomorphic functions 
on $\Mg$ if $g \geq 3$, while the $3g-3$-dimensional orbifold $\Mg$ has an enormous amount of $\smooth$ functions. Hence we would like to find functions which have some geometric and/or arithmetic meanings. 
Hyperbolic geometry is a highly effective tool for studying Riemann surfaces and Teichm\"uller spaces. But, in the context of hyperbolic geometry, the author knows no interesting functions on the moduli space $\Mg$ of compact Riemann surfaces except the injective radius of a compact hyperbolic surface. \par
Other effective tools for studying compact Riemann surfaces are the classical theory of Abelian integrals and its powerful descendant, arithmetic geometry including the Arakelov geometry. Theta functions provide meaningful functions 
on the modulli space $\Mg$. Moreover, 
based on his own Arakelov geometry, Faltings \cite{Fa84} introduced an analytic invariant of a compact Riemann surface. Hain and Reed \cite{HR2004} introduced another invariant by the Hodge line bundle and the intermediate Jacobian on the moduli space 
$\Mg$. Zhang \cite{Zhang}  and the author \cite{Kaw08pre, Kaw2009} 
introduced the same invariant $a_g = \varphi$ independently.
Zhang's approach comes from the Arakelov geometry, while
the author's from the Johnson homomorphism similar to that of Hain and Reed, but based on the classical theory of Abelian integrals. Later de Jong \cite{deJ2013} proved 
that the Faltings invariant, the Hain-Reed invariant and the Kawazumi-Zhang invariant are linearly dependent. Moreover he has clarified various aspects of the invariant 
$a_g$ in papers including \cite{deJ2014, deJ2015, deJ2016}. 
In particular, de Jong and Shokrieh \cite[Theorem 7.1]{JS} gave a complete description of the asymptotic behavior of the invariant $a_g$ along $1$-parameter families of compact Riemann surfaces degenerating into a stable curve.  
D'Hoker and Green \cite{DG2014} found some relation of the genus $2$ invariant 
$a_2$ to high energy physics. D'Hoker, Green, Pioline and Russo \cite{DGPR} proved the function $a_2$ is an eigenfunction with respect to the Laplacian on the 
moduli space of principally polarized Abelian varieties $\mathbb{A}_2$. 
In \S\ref{subsec:jac} we will show this does not hold in the case $g=3$
by using a twisted version of the invariant $a_g$. 
The obstruction for $a_3$ to be an eigenfunction 
is the (extended) first Johnson homomorphism \cite{J80b, Mo93}. 
\par
In this paper, we introduce a twisted version of the invariant $a_g$, 
a map $A$ from the $4$-fold tensor product of the first cohomology group 
of the surface to the complex numbers. One can recover the invariant $a_g$ 
from the map $A$ by using the intersection pairing on the cohomology group, 
and can construct new invariants by a similar way to $a_g$. 
Moreover, from the invariant $A$, we obtain an explicit $(1,1)$-form on the moduli space $\Mg$ representing the first Mumford-Morita-Miller class $e_1 = \kappa_1$
in \S\ref{subsec:e_1}.\par
The invariant $A$ is given as an integral over the $2$-fold product of the surface. 
D'Hoker and Schlotterer \cite{DS2022} generalized our construction to 
integrals over the $n$-fold product of the surface for $n \geq 3$, and 
proved highly nontrivial identities among them. 
It is amazing that their proof of the identities 
is based only on formal natures of the Arakelov Green function $G(-,-)$.
\par
This paper is organized as follows. The twisted version $A$ of the invariant $a_g$ of a compact Riemann surface $C$ is introduced in \S\ref{sec:twisted}.
In \S\ref{sec:tensor} we give two kinds of interpretation of the map $A$  as tensor fields on the moduli space $\Mg$.
First it defines a $(1,1)$-form on $\Mg$. As is proved in \S\ref{sec:e_1}, it is closed and represents the first Mumford-Miller-Morita class $e_1 = \kappa_1$.
In \S\ref{subsec:jac}, using the second interpretation, we discuss how the invariant $a_g$ behaves 
in the cases $g= 2,3$. The proof of Proposition \ref{prop:a3}, which we need for our consideration in the case $g=3$, is given in \S\ref{sec:a3}. 

\par\bigskip\noindent
{\bf Acknowledgements.}\,\, 
First of all the author expresses his sincere thanks and respect 
to Athanase Papadopoulos for his friendship and 
outstanding leadership in world-wide research networks on Teichm\"uller theory.
As well as his own research, many books he has edited illustrate his great insight. The author also thanks M. Uludag 
for giving him a chance to take part in a tribute to 
Athanase Papadopoulos. 
This paper is based on the author's talks at OIST, University of Geneva and 
University of Strasbourg in 2016 and 2017. The slides for OIST and 
Strasbourg are uploaded in \cite{Slides1, Slides2}. The author thanks 
S. Hikami, A. Alekseev and C. Vespa for giving him these opportunities.
In particular, discussions at OIST with Pioline aroused 
the author's interest in the genus $3$ invariant $a_3$. He deeply thanks 
B. Pioline. Furthermore the author thanks R. de Jong for his 
helpful comments on the first version of this paper. 
The present research is supported in part by the grants JSPS KAKENHI 
18KK0071, 18K03283, 19H01784, 20H00115 and 22H01120. 
\par

%\tableofcontents
%\newpage

\section{A twisted invariant of a compact Riemann surface}
\label{sec:twisted}

Throughout this paper, we assume the genus $g$ of a compact Riemann surface $C$ is not smaller than $1$, $g \geq 1$. Then the vector space 
$H' := H^0(C; \Omega^1_C)$ of holomorphic $1$-forms on $C$ does not vanish.
Its complex conjugate $H'' := \overline{H^0(C; \Omega^1_C)}$ is the vector space of anti-holomorphic $1$-forms on $C$. 
The first cohomology group $H := H^1(C; \bC)$ is naturally identified with 
the direct sum $H'\oplus H''$, the vector space of harmonic $1$-forms on the surface $C$. The intersection pairing
\begin{equation}\label{eq:cdot}
(\phi_1, \phi_2) \in H\times H \mapsto 
\phi_1\cdot\phi_2 := \int_C \phi_1\wedge \phi_2 \in \bC,
\end{equation}
is non-degenerate. Hence we identify $H$ with its dual through the pairing.
Under this identification, the subspaces $H'$ and $H''$ are dual to each other. 
\par
For $\phi'_1, \phi'_2 \in H' = H^0(C; \Omega^1_C)$ we can define a positive definite Hermitian pairing 
\begin{equation}\label{eq:pair}
(\phi'_1, \phi'_2) \in H'\otimes H' \mapsto \frac{\sqrt{-1}}2\int_C \phi'_1\wedge\overline{\phi'_2} \in \bC.
\end{equation}
In fact, any holomorphic $1$-form $\phi'$ locally given by $\phi' = f(z)dz$ 
in terms of a local complex coordinate $z = x + \sqrt{-1}y$
satisfies $\frac{\sqrt{-1}}2\phi'\wedge\overline{\phi'} = |f(z)|^2dx\wedge dy$. 
We take an orthonormal basis $\{\psi_i\}^g_{i=1} \subset H'$ 
with respect to the Hermitian pairing
$$
\frac{\sqrt{-1}}2\int_C \psi_i\wedge\overline{\psi_j} = \delta_{ij}, \quad 
(1 \leq i, j \leq g).
$$
Then the element
\begin{equation}\label{eq:omega0}
\om := \sum^g_{i=1}\psi_i\otimes \overline{\psi_i} \in H'\otimes H''
(\subset H\otimes H)
\end{equation}
is independent of the choice of the basis $\{\psi_i\}^g_{i=1}$. 
In fact, the set $\{\psi_i, \frac{\sqrt{-1}}{2}\overline{\psi_i}\}_{1\leq i \leq g}$ is 
a symplectic basis of $H = H^1(C; \bC)$ with respect to the pairing \eqref{eq:cdot}.
Hence the symplectic form on $H$, which depends only on the pairing \eqref{eq:cdot}, equals 
$$
\frac{\sqrt{-1}}{2}\sum^g_{i=1}\psi_i\otimes \overline{\psi_i} 
- \frac{\sqrt{-1}}{2}\sum^g_{i=1}\overline{\psi_i} \otimes \psi_i \in H^1(C; \bC)^{\otimes 2}, 
$$
whose component in $H'\otimes H''$ equals $\frac{\sqrt{-1}}{2}\om$. %\par
In particular, a $(1,1)$-form $B$ on the surface $C$ defined by 
$$
B := \frac{\sqrt{-1}}{2g}\sum^g_{i=1}\psi_i\wedge\overline{\psi_i},
$$
is independent of the choice of the basis $\{\psi_i\}^g_{i=1}$. 
We have $\int_CB = 1$. Moreover the form $B$ is a real volume form on the surface $C$, since the $1$-forms $\psi_i$'s have no common zeroes
because of the Riemann-Roch formula.
%the holomorphic line bundle $\Omega^1_C$ is basepoint free. 
\par
We denote by $\Omega^p(C)$ the vector space of $\bC$-valued $p$-forms on $C$ for $p = 0,1,2$. 
Then the Hodge $*$-operator $*: \Omega^1(C) \to \Omega^1(C)$ depends only on 
the complex structure of the surface $C$. In fact, in terms of a local complex coordinate $z$, we have $*dz = -\sqrt{-1}dz$ and $*d\zb = \sqrt{-1}d\zb$.
Then, since the operator $d*d$ is equivalent to the Laplacian by means of the volume form $B$, we have an exact sequence 
$$
0 \to \bC \to \Omega^0(C) \overset{d*d}\longrightarrow \Omega^2(C) 
\overset{\int_C}\longrightarrow \bC \to 0,
$$
where the left $\bC$ means the constant functions. 
Hence we have a unique linear map
$
\hP: \Omega^2(C) \to \Omega^0(C)
$
such that 
$$
d*d\hP(\Omega) = \Omega - \left(\int_C\Omega\right)B\quad\text{and}\quad
\int_C\hP(\Omega)B = 0
$$
for any $\Omega \in \Omega^2(C)$. It is nothing but the Green operator associated 
with the volume form $B$. The operator $\hP$ is real, $\overline{\hP} = \hP$,
and the induced map
$$
\hP\vert_{\Ker(\int_C)}: \Ker\left(\int_C: \Omega^2(C) \to \bC\right)
\to \Omega^0(C)/\bC
$$
depends only on the complex structure on the surface $C$. 
Moreover we have
\begin{equation}\label{eq:green}
\int_C\hP(\Omega)\Omega' = \int_C\Omega\hP(\Omega') 
= \frac1{4\pi}\int\!\!\!\!\int_{(P_1,P_2)\in C\times C}
(\Omega)_{P_1}(\log G(P_1, P_2))(\Omega')_{P_2}
\end{equation}
for any $\Omega, \Omega' \in \Omega^2(C)$, 
where we denote by $G(-,-)$ the Arakelov Green function.
Here we remark $\hP(B) = 0$. In fact, the function
$\hP(B) $ is constant since $d*d\hP(B) = 0$, and so $\hP(B) = 
\int_C\hP(B)B = 0$. 
\par
Now we define a map $A$ by
\begin{equation}\label{eq:mapA}
A: H'\otimes H'' \otimes H'\otimes H''  \to \bC, \quad
 \phi_1\otimes \overline{\phi_2}\otimes \phi_3\otimes \overline{\phi_4}
\mapsto 
\int_C \phi_1\wedge \overline{\phi_2} \,\hP(\phi_3\wedge
\overline{\phi_4}),
\end{equation}
which is independent of the choice of the basis $\{\psi_i\}^g_{i=1}$, and so 
is a twisted invariant of the compact Riemann surface $C$. 
One can extend it to a map $H\otimes H\otimes H\otimes H \to \bC$ 
by the same formula \eqref{eq:mapA}, but with no additional information.
In order to compute the map $A$ explicitly, it is convenient to 
introduce a complex number $A_{i\jb k \lb}$ defined by 
\begin{equation}\label{eq:Aijkl}
\aligned
A_{i\jb k \lb} & := A(\psi_i\otimes \overline{\psi_j}\otimes \psi_k\otimes \overline{\psi_l}) = 
\int_C \psi_i\wedge \overline{\psi_j} \,\hP(\psi_k\wedge
\overline{\psi_l}) \\
& = \frac1{4\pi}\int\!\!\!\!\int_{(P_1,P_2)\in C\times C}
(\psi_i\wedge \overline{\psi_j})_{P_1}(\log G(P_1, P_2))(\psi_k\wedge
\overline{\psi_l})_{P_2}
\endaligned
\end{equation}
for $1 \leq i, j, k, l \leq g$ by using the orthonormal basis $\{\psi_i\}^g_{i=1}$. Then we have 
\begin{equation}
\aligned
& \overline{A_{i\jb k \lb}} = \int_C \psi_j\wedge \overline{\psi_i} \,\hP(\psi_l\wedge\overline{\psi_k}) = A_{j\ib l \kb}, \\
& A_{i\jb k \lb} = \int_C \hP(\psi_i\wedge \overline{\psi_j}) \, \psi_k\wedge
\overline{\psi_l} = A_{k \lb i\jb}.
\endaligned
\label{eq:elem}
\end{equation}
The latter follows from \eqref{eq:green}.
We can contract the indices $i, k, \dots$ and the ones 
$\jb, \lb, \dots$ to get a new complex number.
For example, the Kawazumi-Zhang invariant\footnote{
The author uses the symbol $a_g = a_g(C)$, while Zhang and others use the symbol $\varphi = \varphi(C)$.} $a_g = a_g(C)$
is obtained by 
$$
a_g = \sum^g_{i,j=1} A_{i\jb j \ib} = A((24)(\om\otimes\om)),
$$
which is a real number by the formula \eqref{eq:elem}. 
Here $(24): H'\otimes H'' \otimes H'\otimes H'' \to H'\otimes H'' \otimes H'\otimes H''$ is the switch map of the 2nd and the 4th components
of the tensor product $H'\otimes H'' \otimes H'\otimes H''$. \par
We can describe the numbers $A_{i\jb k \lb}$ 
and $a_g$ in the diagrams\par\medskip
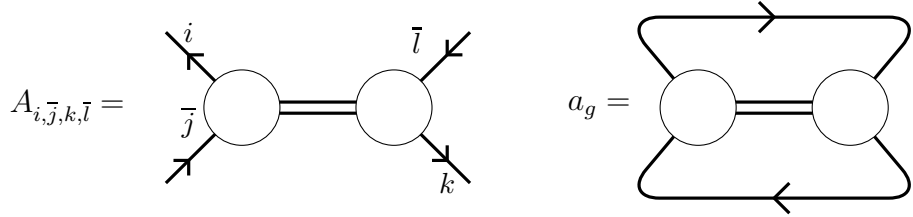
\begin{tikzpicture}
\tikzset{Process/.style={rectangle,  draw,  text centered, text width=37mm, minimum height=1.5cm}};
\draw (-2.3,0) node{\large $A_{i, \overline{j}, k, \overline{l}} = $};
\draw[very thick, double distance=3pt] (0,0) -- (2,0);
\draw[very thick] (-1,1) -- (0,0);
\draw[very thick] (-0.7,0.7) -- (-0.5,0.7);
\draw[very thick] (-0.7,0.7) -- (-0.7,0.5);
\draw (-0.7,1) node{$i$};
\draw[very thick] (-1,-1) -- (0,0);
\draw[very thick] (-0.7,-0.7) -- (-0.9,-0.7);
\draw[very thick] (-0.7,-0.7) -- (-0.7,-0.9);
\draw (-0.7,-0.2) node{$\overline{j}$};
\draw[very thick] (2,0) -- (3,1);
\draw[very thick] (2.7,0.7) -- (2.9,0.7);
\draw[very thick] (2.7,0.7) -- (2.7,0.9);
\draw (2.3,0.9) node{$\overline{l}$};
\draw[very thick] (2,0) -- (3,-1);
\draw[very thick] (2.7,-0.7) -- (2.5,-0.7);
\draw[very thick] (2.7,-0.7) -- (2.7,-0.5);
\draw[fill=white] (0,0)  circle [radius=.5] ;
\draw[fill=white] (2,0)  circle [radius=.5] ;
\draw (2.7,-1) node{$k$};
%\tikzset{Process/.style={rectangle,  draw,  text centered, text width=37mm, minimum height=1.5cm}};
\draw (4.7,0) node{\large $a_g = $};
\draw[very thick, double distance=3pt] (6,0) -- (8,0);
\draw[very thick] [rounded corners=14pt] (8,0) -- (9,1.2) -- (5,1.2) -- (6,0) ;
\draw[very thick] (7,1.2) -- (6.8,1.4);
\draw[very thick] (7,1.2) -- (6.8,1.0);
\draw[very thick] [rounded corners=14pt] (8,0) -- (9,-1.2) -- (5, -1.2) -- (6,0);
\draw[very thick] (7,-1.2) -- (7.2,-1.4);
\draw[very thick] (7,-1.2) -- (7.2,-1.0);
\draw[fill=white] (6,0)  circle [radius=.5] ;
\draw[fill=white] (8,0)  circle [radius=.5] ;
\end{tikzpicture}.
\par\medskip\noindent
Here the circle with indices $i$ and $\jb$ means $\psi_i\wedge\overline{\psi_j}$, 
the double line without arrows $\log G(-,-)$, and a line with an arrow the contraction. For example, since $\frac{\sqrt{-1}}{2g}\hP(\sum^g_{i=1}\psi_i\wedge \overline{\psi_i}) = \hP(B) = 0$, we have 
$\sum_{i, j,l =1}^gA_{i\jb l \lb} = \sum_{j,k,l =1}^gA_{j\jb k \lb} = 0$, 
that is, 
\par
\medskip
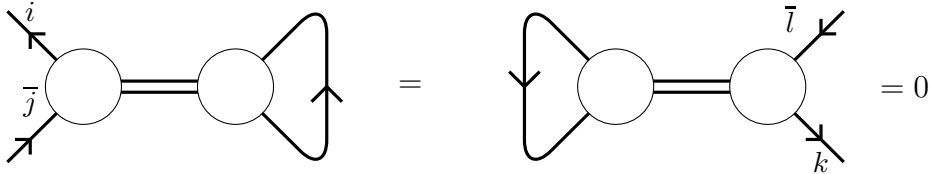
\begin{tikzpicture}
\tikzset{Process/.style={rectangle,  draw,  text centered, text width=37mm, minimum height=1.5cm}};
%\draw (-2.3,0) node{\large $A_{i, \overline{j}, k, \overline{l}} = $};
\draw[very thick, double distance=3pt] (0,0) -- (2,0);
\draw[very thick] (-1,1) -- (0,0);
\draw[very thick] (-0.7,0.7) -- (-0.5,0.7);
\draw[very thick] (-0.7,0.7) -- (-0.7,0.5);
\draw (-0.7,1) node{$i$};
\draw[very thick] (-1,-1) -- (0,0);
\draw[very thick] (-0.7,-0.7) -- (-0.9,-0.7);
\draw[very thick] (-0.7,-0.7) -- (-0.7,-0.9);
\draw (-0.7,-0.2) node{$\overline{j}$};
%\draw[very thick] (2,0) -- (3,1);
%\draw[very thick] (2.7,0.7) -- (2.9,0.7);
%\draw[very thick] (2.7,0.7) -- (2.7,0.9);
%\draw (2.3,0.9) node{$\overline{l}$};
\draw[very thick] [rounded corners=14pt] (2,0) -- (3.2,-1.2) -- (3.2,1.2) -- (2,0);
\draw[very thick] (3.2,0) -- (3,-0.2);
\draw[very thick] (3.2,0) -- (3.4,-0.2);
\draw[fill=white] (0,0)  circle [radius=.5] ;
\draw[fill=white] (2,0)  circle [radius=.5] ;
\draw (4.3,0) node{\large $=$};
\draw[very thick, double distance=3pt] (7,0) -- (9,0);
\draw[very thick] [rounded corners=14pt] (7,0) -- (5.8,-1.2) -- (5.8, 1.2) -- (7,0);
\draw[very thick] (5.8,0) -- (5.6,0.2);
\draw[very thick] (5.8,0) -- (6.0,0.2);
\draw[very thick] (9,0) -- (10,1);
\draw[very thick] (9.7,0.7) -- (9.9,0.7);
\draw[very thick] (9.7,0.7) -- (9.7,0.9);
\draw (9.3,0.9) node{$\overline{l}$};
\draw[very thick] (9,0) -- (10,-1);
\draw[very thick] (9.7,-0.7) -- (9.5,-0.7);
\draw[very thick] (9.7,-0.7) -- (9.7,-0.5);
\draw[fill=white] (7,0)  circle [radius=.5] ;
\draw[fill=white] (9,0)  circle [radius=.5] ;
\draw (9.7,-1) node{$k$};
\draw (10.8,0) node{\large $= 0$};
\end{tikzpicture}\par\medskip\noindent
Similarly we can introduce complex numbers associated with the Riemann surface $C$, i.e., functions on the moduli space $\Mg$ of compact Riemann surfaces.
For example, we can consider a function on $\Mg$ defined by the diagram
\par\medskip
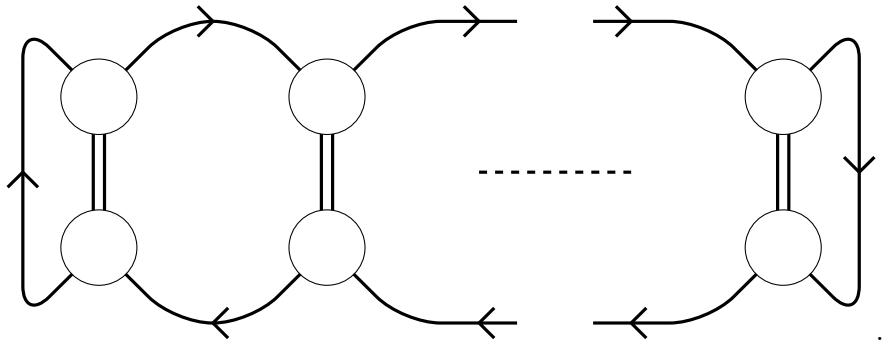
\begin{tikzpicture}
\tikzset{Process/.style={rectangle,  draw,  text centered, text width=37mm, minimum height=1.5cm}};
\draw[very thick, double distance=3pt] (0,0) -- (0,2);
\draw[very thick][rounded corners=14pt] (0,0) -- (-1,-1) -- (-1,3) -- (0,2);
\draw[very thick] (-1,1) -- (-0.8,0.8);
\draw[very thick] (-1,1) -- (-1.2,0.8);
\draw[very thick] [rounded corners=14pt] (0,2) -- (1,3) -- (2,3) -- (3,2);
\draw[very thick] (1.5,3) -- (1.3,2.8);
\draw[very thick] (1.5,3) -- (1.3,3.2);
\draw[very thick] [rounded corners=14pt] (0,0) -- (1,-1) -- (2,-1) -- (3,0) ;
\draw[very thick] (1.5,-1) -- (1.7,-0.8);
\draw[very thick] (1.5,-1) -- (1.7,-1.2);
\draw[fill=white] (0,0)  circle [radius=.5] ;
\draw[fill=white] (0,2)  circle [radius=.5] ;
\draw[very thick, double distance=3pt] (3,0) -- (3,2);
\draw[very thick] [rounded corners=14pt] (3,0) -- (4,-1) -- (5.5, -1);
\draw[very thick] (5,-1) -- (5.2,-0.8);
\draw[very thick] (5,-1) -- (5.2,-1.2);
\draw[very thick] [rounded corners=14pt]  (3,2) -- (4,3) -- (5.5,3);
\draw[very thick] (5,3) -- (4.8,3.2);
\draw[very thick] (5,3) -- (4.8,2.8);
\draw[fill=white] (3,0)  circle [radius=.5] ;
\draw[fill=white] (3,2)  circle [radius=.5] ;
\draw[very thick, dashed] (5,1) -- (7,1) ;
\draw[very thick, double distance=3pt] (9,0) -- (9,2);
\draw[very thick] [rounded corners=14pt] (6.5, -1) -- (8,-1) -- (9,0);
\draw[very thick] (7,-1) -- (7.2,-1.2);
\draw[very thick] (7,-1) -- (7.2,-0.8);
\draw[very thick] [rounded corners=14pt]  (9,2) -- (10,3) -- (10,-1) -- (9,0);
\draw[very thick] (10,1) -- (9.8,1.2);
\draw[very thick] (10,1) -- (10.2,1.2);
\draw[very thick] [rounded corners=14pt] (9,2) -- (8,3) -- (6.5,3) ;
\draw[very thick] (7,3) -- (6.8,3.2);
\draw[very thick] (7,3) -- (6.8,2.8);
\draw[fill=white] (9,0)  circle [radius=.5] ;
\draw[fill=white] (9,2)  circle [radius=.5] ;
\end{tikzpicture}.
\par\medskip
The number $A_{i\jb k \lb}$ is given by an integral over the $2$-fold direct product
$C^{\times 2} = C\times C$. We can generalize this construction to the $3$-fold direct product and more. For example, we can consider
$$
\aligned
&\int_C\hP(\psi_i\wedge\overline{\psi_j})\psi_m\wedge\overline{\psi_n}\, 
\hP(\psi_k\wedge\overline{\psi_l}) 
 = \int_C \psi_i\wedge\overline{\psi_j}\, \hP((\psi_m\wedge\overline{\psi_n})\, 
\hP(\psi_k\wedge\overline{\psi_l}))\\
& = \frac1{(4\pi)^2}\int\!\!\!\!\int\!\!\!\!\int_{(P_1, P_2, P_3) \in C\times C\times C}
(\psi_i\wedge \overline{\psi_j})_{P_1}(\log G(P_1, P_2))(\psi_m\wedge
\overline{\psi_n})_{P_2}(\log G(P_2, P_3))(\psi_k\wedge \overline{\psi_l})_{P_3}
\endaligned
$$
for $1 \leq i, j, k, l, m, n \leq g$. This can be described in the diagram\par
\medskip
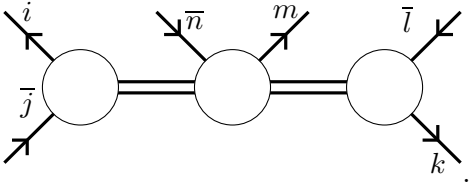
\begin{tikzpicture}
\tikzset{Process/.style={rectangle,  draw,  text centered, text width=37mm, minimum height=1.5cm}};
\draw[very thick, double distance=3pt] (0,0) -- (2,0);
\draw[very thick, double distance=3pt] (2,0) -- (4,0);
\draw[very thick] (-1,1) -- (0,0);
\draw[very thick] (-0.7,0.7) -- (-0.5,0.7);
\draw[very thick] (-0.7,0.7) -- (-0.7,0.5);
\draw (-0.7,1) node{$i$};
\draw[very thick] (-1,-1) -- (0,0);
\draw[very thick] (-0.7,-0.7) -- (-0.9,-0.7);
\draw[very thick] (-0.7,-0.7) -- (-0.7,-0.9);
\draw (-0.7,-0.2) node{$\overline{j}$};
\draw[very thick] (2,0) -- (3,1);
\draw[very thick] (2.7,0.7) -- (2.5,0.7);
\draw[very thick] (2.7,0.7) -- (2.7,0.5);
\draw[very thick] (2,0) -- (1,1);
\draw[very thick] (1.3,0.7) -- (1.3,0.9);
\draw[very thick] (1.3,0.7) -- (1.1,0.7);
\draw (1.5,0.9) node{$\overline{n}$};
\draw (2.7,1) node{$m$};
\draw[very thick] (4,0) -- (5,1);
\draw[very thick] (4.7,0.7) -- (4.9,0.7);
\draw[very thick] (4.7,0.7) -- (4.7,0.9);
\draw[very thick] (4,0) -- (5,-1);
\draw[very thick] (4.7,-0.7) -- (4.5,-0.7);
\draw[very thick] (4.7,-0.7) -- (4.7,-0.5);
\draw (4.3,0.9) node{$\overline{l}$};
\draw (4.7,-1) node{$k$};
\draw[fill=white] (0,0)  circle [radius=.5] ;
\draw[fill=white] (2,0)  circle [radius=.5] ;
\draw[fill=white] (4,0)  circle [radius=.5] ;
\end{tikzpicture}.\par\medskip\noindent
%After the slides \cite{Slides1, Slides2}, 
D'Hoker and Schlotterer \cite{DS2022} gave a systematic generalization of $A_{i\jb k \lb}$ where they constructed modular graph tensors using integrals on the direct product $C^{\times n}$ for any $n \geq 1$, and proved highly nontrivial identities among them. It is amazing that their proof of the identities 
is based only on formal natures of the Arakelov Green function $G(-,-)$. \par

\section{Tensor fields on the moduli space $\Mg$}
\label{sec:tensor}
When the Riemann surface $C$ runs over the moduli space $\Mg$ of compact 
Riemann surfaces of genus $g$, 
the transpose ${}^tA$ of the linear map $A$ in \eqref{eq:mapA} can be regarded
as some tensor fields on the moduli space $\Mg$. 
Now recall the vector spaces $H'$ and $H''$ are dual to each other
through the pairing \eqref{eq:cdot}. In particular, 
the basis $\{\psi_i\}^g_{i=1}$ and $\{\frac{\sqrt{-1}}2\overline{\psi_i}\}^g_{i=1}$ are 
dual to each other. Then the transpose ${}^tA: \bC \to H''\otimes H'\otimes H'' \otimes H'$ is determined by its own value at $1 \in \bC$
\begin{equation}\label{eq:A*}
{}^tA(1) = \frac1{16}\sum_{i,j,k,l = 1}^g A_{i\jb k \lb} 
\overline{\psi_i}\otimes \psi_j\otimes \overline{\psi_k}\otimes \psi_l
\in H''\otimes H'\otimes H'' \otimes H'.
\end{equation}
When the surface $C$ runs over the moduli space $\Mg$, the collections of $H'$ and $H''$ are holomorphic and anti-holomorphic vector bundles, respectively.
Then ${}^tA(1)$ is a $\smooth$ section of a tensor product of these 
vector bundles. \par
\subsection{The first Mumford-Morita-Miller class $e_1 = \kappa_1$}
\label{subsec:e_1}
Here we recall the holomorphic cotangent space $T^*_{[C]}\Mg$ of the moduli space $\Mg$ 
at the equivalence class $[C]$ of the Riemann surface $C$ is canonically 
isomorphic to the space of holomorphic quadratic differentials
$$
T^*_{[C]}\Mg = H^0(C; (\Omega^1_C)^{\otimes 2}).
$$
For example, we may regard $\psi_l\psi_j = \psi_j\psi_l \in T^*_{[C]}\Mg$ 
and $\overline{\psi_i}\overline{\psi_k} = \overline{\psi_k}\overline{\psi_i} \in \overline{T^*_{[C]}\Mg}$. Hence we can define a $(1,1)$-form $e_1^{\hP}$ 
on the moduli space $\Mg$ by\footnote{In \cite{Slides1, Slides2} 
we denote $e_1^{\hP}$ by $e_1^D$, whose
coefficients was $48\sqrt{-1}$. It is not true. This mistake comes from 
computation of the intersection number $Y_i\cdot \overline{Y_i}$.} 
$$
e_1^{\hP}\vert_{[C]} := -3\sqrt{-1}\sum_{i,j,k,l =1}^g
(\psi_j\psi_l A_{i\jb k \lb} \overline{\psi_k}\overline{\psi_i} 
- \psi_j\psi_l A_{i\jb k \ib} \overline{\psi_k}\overline{\psi_l}).
$$
In fact, from \eqref{eq:A*}, the RHS is independent of the choice of the basis $\{\psi_i\}^g_{i=1}$. 
Then we have 
\begin{theorem}\label{thm:e_1}
The $(1,1)$-form $e_1^{\hP}$ is closed on the moduli space $\Mg$, 
and represents %$\dfrac1{12}$ times 
the first Mumford-Morita-Miller class $e_1 \in H^2(\Mg; \bC)$. 
\end{theorem}
The theorem was stated implicitly in the equation (5.5) of the preprint 
\cite{Kaw08pre} without proof. We will give its proof in the next section \S\ref{sec:e_1}.\par
Very recently, de Jong and van der Lugt \cite{JL, Lugt} introduced a subalgebra 
$\mathcal{R}(\Mg)$ of the de Rham complex $\Omega^*(\Mg)$, called 
{\it the ring of tautological forms} on $\Mg$.
It consists entirely of closed forms, and its image in the cohomology algebra 
$H^*(\Mg; \bC)$ is generated by all the Mumford-Morita-Miller classes, i.e., 
equals the image of the tautological ring of $\Mg$.
It is remarkable that the ring $\mathcal{R}(\Mg)$ is finite-dimensional,
which suggests us that study of differential forms on $\Mg$ could help us to understand the tautological ring of $\Mg$. 
As was shown in \cite{JL, Lugt}, the degree $2$ part of the ring $\mathcal{R}(\Mg)$ is closely related to the invariant $a_g$. So it would be very interesting if we could find further relations of the ring $\mathcal{R}(\Mg)$ with 
the complex numbers $A_{i\jb k \lb}$'s or generalizations by D'Hoker and Schlotterer \cite{DS2022}.
\par
\subsection{The Jacobian map $\Jac: \Mg \to \Ag$}\label{subsec:jac}
The transpose ${}^tA(1)$ in \eqref{eq:A*} has another interpretation: 
it can be regarded as a $\smooth$ section of the pullback $\Jac^*(T^*\Ag\otimes 
\overline{T^*\Ag})$ over the moduli space $\Mg$. 
Here $\Ag$ is the moduli space of principally polarized Abelian varieties 
of genus $g$, and $\Jac: \Mg \to \Ag$, $[C] \mapsto [\Jac(C)]$, the Jacobian map. The moduli space $\Ag$ is the quotient space of the Siegel upper half space $\Hg$ by the Siegel modular group $Sp_{2g}(\bZ)$. The space $\Hg = Sp_{2g}(\bR)/U_g$ is the space of all complex structures $J$ on the standard symplectic vector space $\bR^{2g}$ with positivity condition. 
The tautological vector bundle $E_g := \coprod_{J \in \Hg}(\bR^{2g}, J)$ 
over $\Hg$ is an $Sp_{2g}(\bR)$-equivariant holomorphic vector bundle with canonical Hermitian metric, and so descends to the moduli space $\Ag$. 
Through a canonical isomorphism $T\Ag = \mathrm{Sym}^2(E_g) \subset {E_g}^{\otimes 2}$, one can define a K\"ahler metric on the moduli space $\Ag$. 
See, for example, \cite[Chapter II]{MuTata}. In particular, we may consider 
the Laplacian $\Delta$ on the moduli space $\Ag$. 
The pullback $\Jac^*E^*_g$ of the Hermitian vector bundle $E_g$ under the Jacobian map $\Jac$ at the equivalence class $[C] \in \Mg$ is isometric to 
the holomorphic $1$-forms $H' = H^0(C; \Omega^1_C)$ with metric given in  \eqref{eq:pair}. Hence the transpose ${}^tA(1)$ in \eqref{eq:A*} 
can be regarded as a $\smooth$ section of the pullback 
$\Jac^*(T^*\Ag\otimes \overline{T^*\Ag}))$. \par

By Rauch's formula \cite{R}, the differential $d\Jac$ of the Jacobian map 
is given by the commutative diagram
\begin{equation}\label{eq:Rauch}
\begin{CD}
T^*_{[C]}\Mg @<{(d\Jac)^*}<< T^*_{[\Jac(C)]}\Ag\\
@| @|\\
H^0(C; (\Omega^1_C)^{\otimes 2}) @<{\text{multiplication}}<< \mathrm{Sym}^2 H^0(C; \Omega^1_C).
\end{CD}
\end{equation}
In particular, from M. Noether's theorem, the differential $d\Jac_{[C]}$ is injective if and only if
$g=2$ or ($g \geq 3$ and $C$ is non-hyperelliptic). On the other hand, 
$\dim \Mg = \dim \Ag$ if and only if $g = 2$ or $3$. Hence the Laplacian $\Delta$ 
on $\Ag$ acts on functions on $\mathbb{M}_2 \overset{\text{open}}\subset 
\mathbb{A}_2$ and $\mathbb{M}_3 \setminus \mathbb{H}_3\overset{\text{open}}\subset \mathbb{A}_3$. Here we denote by $\mathbb{H}_g \subset \Mg$ the locus of hyperelliptic Riemann surfaces of genus $g \geq 2$. 

%\par
%Now recall the Laplacian $\Delta$ on $\Ag$ acts on functions on $\Mg$ if $g = 2$ or $3$.
D'Hoker, Green, Pioline and Russo \cite{DGPR} proved the following remarkable theorem.
\begin{theorem}[\cite{DGPR}(4.23)] $\Delta a_2 = 5a_2$. In particular, the invariant 
$a_2$ is an eigenfunction of the Laplacian $\Delta$. 
\end{theorem}
Based on the theorem, Pioline \cite{P2016} gave an explicit formula for $a_2$ 
in terms of theta functions. Later, for any $g\geq 1$, Wilms \cite{Wilms} gave an explicit formula for $a_g$ in terms of a function introduced by de Jong as well as theta functions. But the methods of Wilms \cite{Wilms} seem to be different from 
those discussed in this paper.\par
There remains the case $g = 3$. In this paper we prove 
\begin{theorem}\label{thm:g3} In the case $g=3$, the invariant $a_3$ is not an eigenfunction of the Laplacian $\Delta$. 
\end{theorem}
The obstruction for $a_3$ to be an eigenfunction 
is the first Johnson homomorphism \cite{J80b, Mo93}, 
which is represented by a twisted $1$-form on $\Mg$ as follows. 
We denote by $U = U(C)$ the kernel of the contraction map 
$\Lambda^3H \to H$, $\phi_1\wedge \phi_2\wedge \phi_3\mapsto 
(\phi_2\cdot\phi_3)\phi_1
+ (\phi_3\cdot\phi_1)\phi_2
+ (\phi_1\cdot\phi_2)\phi_3$. 
It should be remarked that $U$ vanishes if and only if $g \leq 2$. %\par
For $\phi_i = \phi'_i + \phi''_i \in H = H'\oplus H''$, $1 \leq i \leq 3$, we define 
$Q(\phi_1, \phi_2, \phi_3) \in \smooth(C; (T^*C)^{\otimes 2})$ by 
\begin{equation}\label{eq:Q}
Q(\phi_1, \phi_2, \phi_3) := -\sqrt{-1}\phi'_1\pa\hP(\phi_2\wedge\phi_3) 
-\sqrt{-1}\phi'_2\pa\hP(\phi_3\wedge\phi_1) 
-\sqrt{-1}\phi'_3\pa\hP(\phi_1\wedge\phi_2).
\end{equation}
The tensor $Q(\phi_1, \phi_2, \phi_3)$ is alternating in 
$\phi_1, \phi_2, \phi_3$, but 
is not holomorphic in general. 
Since $d*d = -2\sqrt{-1}\pab\pa$, we have 
$$
\aligned
& 2\pab Q\\
= &
\phi'_1(\phi_2\wedge\phi_3 - (\int_C\phi_2\wedge\phi_3)B)
+\phi'_2(\phi_3\wedge\phi_1 - (\int_C\phi_3\wedge\phi_1)B)\\
&\hskip 40mm
+\phi'_3(\phi_1\wedge\phi_2 - (\int_C\phi_1\wedge\phi_2)B)\\
= & \phi'_1\phi'_2\phi''_3 - \phi'_1\phi'_3\phi''_2
+ \phi'_2\phi'_3\phi''_1 - \phi'_2\phi'_1\phi''_3
+ \phi'_3\phi'_1\phi''_2 - \phi'_3\phi'_2\phi''_1\\
&\hskip 40mm - (\phi'_1(\phi_2\cdot\phi_3)
+\phi'_2(\phi_3\cdot\phi_1)
+\phi'_3(\phi_1\cdot\phi_2))B\\
=& - (\phi'_1(\phi_2\cdot\phi_3)
+\phi'_2(\phi_3\cdot\phi_1)
+\phi'_3(\phi_1\cdot\phi_2))B.
\endaligned
$$
This implies the restriction of $Q$ to the subspace $U$ has its value 
in the space $H^0(C; (\Omega^1_C)^{\otimes 2})$ of holomorphic quadratic differentials
$$
Q: U = U(C) \to H^0(C; (\Omega^1_C)^{\otimes 2}) = T^*_{[C]}\Mg. 
$$
We remark the map $Q$ is equivariant under the action of the holomorphic automorphism group $\mathrm{Aut}(C)$. 
In particular, if the Riemann surface $C$ is hyperelliptic, we have 
\begin{equation}\label{eq:hypell}
Q(U(C)) \subset (\text{the $(-1)$-eigenspace of the hyperelliptic involution in $T^*_{[C]}\Mg$}),
\end{equation}
since the involution acts on the space $U(C)$ by $-1_{U(C)}$. \par
We may identify the map $Q$ with its transpose\footnote{For simplicity, we use the same symbol $Q$, not ${}^tQ$.}
$Q: T_{[C]}\Mg \to (U(C))^*$, which 
can be regarded as a twisted $(1,0)$ form on the moduli space $\Mg$. 
Here the pairing 
\begin{equation}\label{eq:M^0}
\aligned
\langle-, -\rangle: \Lambda^3H\times \Lambda^3H &\to \bC, \\
(\phi^{(1)}_1\wedge \phi^{(1)}_2\wedge \phi^{(1)}_3, 
\phi^{(2)}_1\wedge \phi^{(2)}_2\wedge \phi^{(2)}_3) &\mapsto 
\sum_{\sigma, \tau \in \mathfrak{S}_3}(\phi^{(1)}_{\sigma(1)}\cdot \phi^{(2)}_{\tau(1)})(\phi^{(1)}_{\sigma(2)}\cdot \phi^{(2)}_{\tau(2)})(\phi^{(1)}_{\sigma(3)}\cdot \phi^{(2)}_{\tau(3)})
\endaligned
\end{equation}
is non-degenerate on $U = U(C)$, and so we identify $U(C)$ with its dual $(U(C))^*$ through the pairing \eqref{eq:M^0}. 
We denote by the same symbol $U$ the flat vector bundle $\coprod_{[C]\in \Mg}U(C)$ over the moduli space $\Mg$. 
The tensor $Q$ equals the $(1,0)$-part of the first variation of the harmonic volumes introduced by Harris \cite[Theorem 5.8]{Harris}. In particular, the sum $Q + \overline{Q}$ is a twisted closed $1$-form with values in the flat vector bundle 
$U$. As was proved in \cite[\S8]{Kaw2006}, the cohomology class of 
$-(Q+\overline{Q})$ equals the extended first Johnson homomorphism
$\tilde k$ \cite{Mo93}
$$
\tilde k = -[Q+\overline{Q}] \in H^1(\Mg; U).
$$
The cohomology group $H^1(\Mg; U)$ is canonically isomorphic to 
the first group cohomology of the mapping class group of a genus $g$ Riemann 
surface with values in the vector space $U = U(C)$. 
The restriction of $\tilde k$ to the Torelli group equals the first Johnson 
homomorphism \cite{J80b}. 
As was proved by Morita \cite{Mo96}, the first 
Mumford-Morita-Miller class $e_1$ equals the contraction $\langle\tilde k, \tilde k\rangle \in H^2(\Mg; \bC)$ of the cup product $\tilde k \cup \tilde k \in H^2(\Mg; U^{\otimes 2})$ up to non-zero constant factor. By the same recipe using the pairing in \eqref{eq:M^0}, we obtain a $(1,1)$-form 
$e_1^J$ on $\Mg$ representing the class $e_1$, which equals 
$\langle Q+\overline{Q}, Q+\overline{Q}\rangle$ up to non-zero constant factor. 
From \cite[Theorem 6.1]{Kaw08pre}, we have 
\begin{equation}\label{eq:Thm6.1}
\frac1{12}e_1^{\hP} =  \frac{-2\sqrt{-1}g}{2(2g+1)}(\pa\pab a_g) + \frac1{12}e_1^J.
\end{equation}
An alternative proof of the formula \eqref{eq:Thm6.1} using the Deligne pairing 
is given by de Jong \cite{deJ2016}. 
If $g=2$, then $U = 0$, so that $e_1^J = 0$. But $e_1^J \neq 0$ for $g \geq 3$.
\par
In \cite[Appendices B and C]{DGPR}, D'Hoker, Green, Pioline and Russo gave suggestive computations for any genus $g \geq 2$. In the equations 
(B.4) and (C.1), they divide 
the exact $(1,1)$ form $\pab\pa a_g$ into two parts $\psi^1_A+\psi_B$ and 
$\psi^2_A+\psi_C$. Their formulae (C.3) and (C.5) say that $\psi^1_A+\psi_B$ and $\psi^2_A+\psi_C$ equal $e_1^{\hP}$ and $e_1^J$ up to non-zero constant factor, respectively. These formulae are equivalent to our equation \eqref{eq:Thm6.1}. 
The formula (C.9) says that $\psi^2_A+\psi_C$ vanishes if $g=2$. 
In order to understand the formula (C.4) 
$$
\Delta a_g\Bigr\vert_{\psi^1_A+\psi_B} = (2g+1)a_g,
$$
we need to introduce a lift of the differential form $e_1^{\hP}$ to 
an element $\widehat{E_1}$ of 
$\smooth(\Mg; \Jac^*(T^*\Ag\otimes \overline{T^*\Ag}))$
defined by
$$
\widehat{E_1}\vert_{[C]} := -3\sqrt{-1}\sum_{i,j,k,l =1}^g
((\psi_j\bullet \psi_l) A_{i\jb k \lb} (\overline{\psi_k}\bullet \overline{\psi_i})
- (\psi_j\bullet \psi_l) A_{i\jb k \ib} (\overline{\psi_k}\bullet \overline{\psi_l})).
$$
Here $\psi_j\bullet \psi_l := \frac12(\psi_j\otimes \psi_l + \psi_l\otimes \psi_j) 
\in \mathrm{Sym}^2H'$ and so on.
Now we recall the Laplacian $\Delta$ on any K\"ahler manifold $X$ is 
given by 
$$
\Delta f = 2\sqrt{-1}\Lambda \pa\pab f
$$
for any $f \in \smooth(X; \bC)$. Here $\Lambda: T^*X \otimes \overline{T^*X}
\to X\times \bC$ is the K\"ahler contraction. In our case $X = \Ag$, we have 
$$
\Lambda \widehat{E_1}\vert_{[C]}
= 6\left(\sum_{j,l =1}^g A_{l\jb j \lb} + \sum_{j,l =1}^gA_{j\jb l \lb} 
- g\sum_{i,j =1}^gA_{i\jb j \ib}
- \sum_{i,l =1}^gA_{i\lb l \ib}\right) = -6ga_g,
$$
since $\int_C\psi_i\wedge\overline{\psi_j} = -2\sqrt{-1}\delta_{ij}$ and 
$\sum_{j,k,l =1}^gA_{j\jb k \lb} = 0$. Applying the K\"ahler contraction
to the equation \eqref{eq:Thm6.1}, we obtain 
$
-\frac{g}{2}a_g = -\frac{g}{2(2g+1)}\Delta a_g + \frac1{12}\Lambda e_1^J
$, and so 
\begin{equation}\label{eq:2g+1}
\Delta a_g = (2g+1)a_g + \frac{g}{6}\Lambda e_1^J,
\end{equation}
if $g = 2,3$. In the case $g = 3$, 
the K\"ahler contraction $\vert \Lambda e_1^J\vert$ with respect to the metric of $\mathbb{A}_3$ satisfies the following. 
\begin{proposition}\label{prop:a3} There is a hyperelliptic Riemann surface $C_0$ 
of genus $3$ 
%with no nontrivial holomorphic automorphism other than the hyperelliptic involution 
such that 
$$
\vert \Lambda e_1^J\vert \longrightarrow + \infty
$$
as a non-hyperelliptic Riemann surface $C$ of genus $3$ goes to the surface $C_0$. 
\end{proposition}
\noindent
The proof will be given in \S\ref{sec:a3}, and is based on a theorem of Harris \cite[Theorem 6.5]{Harris}.\par
Now the invariant $a_3$ is $\smooth$ near $[C_0] \in \mathbb{M}_3$, so is locally bounded near $[C_0]$. Proposition \ref{prop:a3} together with 
\eqref{eq:2g+1} implies that 
$\Delta a_3$ is not bounded near $[C_0]$. Hence the function 
$a_3$ is not an eigenfunction of the Laplacian $\Delta$. This completes the proof of Theorem \ref{thm:g3} modulo that of Proposition \ref{prop:a3}. 
\qed
\par
\bigskip
Since the Jacobian map $\Jac\vert_{\mathbb{H}_g}: 
\mathbb{H}_g \to \Ag$ restricted to the hyperelliptic locus $\mathbb{H}_g$ is an embedding, the K\"ahler metric on $\Ag$ 
induces a K\"ahler metric on $\mathbb{H}_g$ and its Laplacian $\Delta_{\mathbb{H}_g}$. Since the $(1,1)$-form $e_1^J$ vanishes along $\mathbb{H}_g$ from \eqref{eq:hypell}, one can clarify $\Delta_{\mathbb{H}_g}(a_g\vert_{\mathbb{H}_g})$ if the K\"ahler contraction $\Lambda$ on $\mathbb{H}_g$ can be computed. \par

\section{A proof of Theorem \ref{thm:e_1}}
\label{sec:e_1}

In this section we prove Theorem \ref{thm:e_1} based on computations 
given in \cite{Kaw08pre}. There is introduced a closed $(1,1)$-form $E_1^D$ on the moduli space $\Mg$ representing $\frac1{12}e_1$ in the equation (5.3).
We will prove the equality $E_1^D = \frac{1}{12}e_1^{\hP}$ to obtain the theorem. 
\par
Our starting point is Lemma 5.5 in \cite{Kaw08pre} which says that 
\begin{equation}\label{eq:Lem5.5}
E^D_1(\la, \overline{\mu}) = 4M\int_C\cH(\omega'_{(1)}\la)\wedge
\omega'_{(1)}\hP d*(\omega''_{(1)}\overline{\ell^\mu} - \overline{\ell^\mu}\omega''_{(1)})
\end{equation}
for any Beltrami differentials $\la, \mu \in \smooth(C; TC\otimes \overline{T^*C})$. \par
First of all we need to explain the notation in the equation \eqref{eq:Lem5.5}: 
$M$, $\omega'_{(1)}$, $\omega''_{(1)}$, $\cH$ and $\ell^\mu$. 
The first complex homology group $H_1(C; \bC)$ admits 
the intersection pairing which we denote by $\cdot: H_1(C; \bC)\times H_1(C; \bC)
\to \bC$. 
The integrand of the RHS is an $H_1(C; \bC)^{\otimes 4}$-valued $(1,1)$-form.
The symbol $M: H_1(C; \bC)^{\otimes 4}\to \bC$ is given by 
$M(Z_1Z_2Z_3Z_4):= (Z_2\cdot Z_3)(Z_4\cdot Z_1)$ for any $Z_i \in H_1(C; \bC)$.
Here and throughout this section we omit the symbol $\otimes$. %\par
The first complex cohomology group $H = H^1(C; \bC)$, 
is embedded into the space of $1$-forms $\Omega^1(C)$ as the subspace of harmonic $1$-forms.
An $H_1(C; \bC)$-valued real harmonic $1$-form $\omega_{(1)} \in \Omega^1(C)\otimes H_1(C; \bC)$ is defined to be the embedding $\in \Hom(H^1(C; \bC), \Omega^1(C)) = \Omega^1(C)\otimes H_1(C; \bC)$. We denote by $\omega'_{(1)}$ and $\omega''_{(1)}$ its $(1,0)$- and $(0,1)$- parts, respectively. If homology classes $Y_i \in H$, $1 \leq i \leq g$, satisfy 
$$
\int_{Y_i}\psi_j = \delta_{i,j}, \quad\text{and}\quad
\int_{Y_i}\overline{\psi_j} = 0
$$
for $1\leq i, j \leq g$, in other words, $\{Y_i, -2\sqrt{-1}\overline{Y_i}\}_{1 \leq i\leq g}$ is the dual basis of the symplectic basis $\{\psi_i, \frac{\sqrt{-1}}2\overline{\psi_i}\}_{1 \leq i\leq g}$, then we have 
\begin{equation}\label{eq:omega1}
\omega_{(1)} = \omega'_{(1)} + \omega''_{(1)}, \quad
\omega'_{(1)} = \sum^g_{i=1}\psi_iY_i, \quad\text{and}\quad
\omega''_{(1)} = \sum^g_{i=1}\overline{\psi_i}\overline{Y_i}.
\end{equation}
Here we remark $Y_i\cdot \overline{Y_j} = \frac{\sqrt{-1}}2\delta_{ij}$. 
The Hodge decomposition on the $1$-forms on $C$ is given by 
$$
\varphi = \cH\varphi + *d\hP d\varphi + d\hP d*\varphi
$$
for any $\varphi \in \Omega^1(C)$. Here $\cH$ is the harmonic projection, which is 
written by 
$$
\cH(\varphi) = \frac{\sqrt{-1}}2\sum^g_{j=1}\left(\left(\int_C\varphi\wedge\overline{\psi_j}\right)\psi_j
-  \left(\int_C\varphi\wedge\psi_j\right)\overline{\psi_j}\right).
$$
For example, $\omega'_{(1)}\lambda$ can be regarded as an $H_1(C; \bC)$-valued $(0,1)$-form by contracting $T^*C$ in $\omega'_{(1)}$ and $TC$ in $\lambda$, 
so that we have 
\begin{equation}\label{eq:Hol}
\cH(\omega'_{(1)}\lambda) =\sum^g_{l=1}\cH(\psi_l\la)Y_l
= \frac{\sqrt{-1}}2\sum^g_{l,j=1}\left(\int_C\psi_j\psi_l\la\right)\overline{\psi_j}Y_l,
\end{equation}
since $\psi_l\la\wedge \psi_j = - \psi_j \wedge \psi_l\la = - \psi_j \psi_l\la$. 
On the other hand, we have 
\begin{equation}\label{eq:Pd*}
\varphi'' - \cH(\varphi'') = *d\hP d\varphi'' + d\hP d*\varphi'' = 2\pab\hP d*\varphi''
\end{equation}
for any $(0,1)$-current $\varphi''$. \par
Finally $\ell^\mu \in \Omega^0(C)\otimes H_1(C; \bC)$ is defined by 
$$
\ell^\mu := 2\hP d*(\omega'_{(1)}\mu).
$$
From \eqref{eq:Pd*} we have 
$$
\aligned
& d*(\omega'_{(1)}\ell^\mu - \ell^\mu\omega'_{(1)})
= -\sqrt{-1}\pab(\omega'_{(1)}\ell^\mu - \ell^\mu\omega'_{(1)})
= \sqrt{-1}(\omega'_{(1)}\wedge \pab\ell^\mu - (\pab\ell^\mu)\wedge \omega'_{(1)})\\
&= 2\sqrt{-1}\omega'_{(1)}\wedge \pab\hP d*(\omega'_{(1)}\mu)
+ 2\sqrt{-1}\pab\hP d*(\omega'_{(1)}\mu)\wedge \omega'_{(1)}\\
&= \sqrt{-1}\omega'_{(1)}\wedge(\omega'_{(1)}\mu - \cH(\omega'_{(1)}\mu))
+ \sqrt{-1}(\omega'_{(1)}\mu - \cH(\omega'_{(1)}\mu))\wedge \omega'_{(1)}\\
&= - \sqrt{-1}\omega'_{(1)}\wedge\cH(\omega'_{(1)}\mu)
- \sqrt{-1}\cH(\omega'_{(1)}\mu)\wedge \omega'_{(1)}.
\endaligned
$$
The last equality follows from 
$$
\aligned
&\omega'_{(1)}\wedge \omega'_{(1)}\mu + \omega'_{(1)}\mu\wedge \omega'_{(1)}
= \sum^g_{i,j=1}(\psi_i\wedge\psi_j\mu + \psi_i\mu\wedge\psi_j)Y_iY_j\\
& = \sum^g_{i,j=1}(\psi_i\psi_j - \psi_j\psi_i)\mu Y_iY_j 
= 0.
\endaligned
$$
Hence, from \eqref{eq:Lem5.5}, we have 
$$
E^D_1(\la, \overline{\mu}) = 4\sqrt{-1}M\int_C\cH(\omega'_{(1)}\la)\wedge
\omega'_{(1)}\hP\left(
\cH(\omega''_{(1)}\overline{\mu})\wedge \omega''_{(1)}
+ \omega''_{(1)}\wedge\cH(\omega''_{(1)}\overline{\mu})
\right),
$$
which is exactly the equation (5.5) in \cite{Kaw08pre}. \par
From \eqref{eq:Hol} we have 
$$
\aligned
& 4\sqrt{-1}M\int_C\cH(\omega'_{(1)}\la)\wedge
\omega'_{(1)}\hP(
\cH(\omega''_{(1)}\overline{\mu})\wedge \omega''_{(1)})\\
&= \sqrt{-1}\sum_{j,l,i,k,m,n=1}^g\left(\int_C\psi_j\psi_l\la\right)
\left(\int_C\overline{\psi_k}\overline{\psi_i}\overline{\mu}\right)
\int_C\overline{\psi_j}\wedge\psi_m\hP(\psi_k\wedge\overline{\psi_n})
Z(Y_lY_m\overline{Y_i}\overline{Y_n})\\
&= -\frac{\sqrt{-1}}4\sum_{j,l,i,k,m,n=1}^g\left(\int_C\psi_j\psi_l\la\right)
\left(\int_C\overline{\psi_k}\overline{\psi_i}\overline{\mu}\right)
A_{m\jb k \overline{n}}
%\int_C\overline{\psi_j}\wedge\psi_m\hP(\psi_k\wedge\overline{\psi_n})
\delta_{mi}\delta_{nl}\\
&= -\frac{\sqrt{-1}}4\sum_{j,l,i,k=1}^g\left(\int_C\psi_j\psi_l\la\right)
\left(\int_C\overline{\psi_k}\overline{\psi_i}\overline{\mu}\right)
A_{i\jb k \lb}.
\endaligned
$$
On the other hand, we have 
$$
\aligned
& 4\sqrt{-1}M\int_C\cH(\omega'_{(1)}\la)\wedge
\omega'_{(1)}\hP(\omega''_{(1)}\wedge\cH(\omega''_{(1)}\overline{\mu}))\\
& = \sqrt{-1}\sum_{j,l,i,k,m,n=1}^g\left(\int_C\psi_j\psi_l\la\right)
\left(\int_C\overline{\psi_k}\overline{\psi_i}\overline{\mu}\right)
\int_C\overline{\psi_j}\wedge\psi_m\hP(\overline{\psi_n}\wedge\psi_k)
Z(Y_lY_m\overline{Y_n}\overline{Y_i})\\
& = \frac{\sqrt{-1}}4\sum_{j,l,i,k,m,n=1}^g\left(\int_C\psi_j\psi_l\la\right)
\left(\int_C\overline{\psi_k}\overline{\psi_i}\overline{\mu}\right)
A_{m\jb k \overline{n}}
%\int_C\overline{\psi_j}\wedge\psi_m\hP(\overline{\psi_n}\wedge\psi_k)
\delta_{mn}\delta_{il}\\
& = \frac{\sqrt{-1}}4\sum_{j,l,k,m=1}^g\left(\int_C\psi_j\psi_l\la\right)
\left(\int_C\overline{\psi_k}\overline{\psi_l}\overline{\mu}\right)
A_{m\jb k \overline{m}}
= \frac{\sqrt{-1}}4\sum_{i,j,l,k=1}^g\left(\int_C\psi_j\psi_l\la\right)
\left(\int_C\overline{\psi_k}\overline{\psi_l}\overline{\mu}\right)
A_{i\jb k \overline{i}}
\endaligned
$$
Consequently we obtain
$$
E^D_1 = -\frac{\sqrt{-1}}4\sum_{i,j,k,l=1}^g\left(\psi_i\psi_lA_{i\jb k \lb}
\overline{\psi_k}\overline{\psi_i} - \psi_j\psi_lA_{i\jb k \overline{i}}
\overline{\psi_k}\overline{\psi_l}
\right) = \frac1{12}e_1^{\hP},
$$
as was to be shown. \qed
\par
\bigskip
The present proof is based on direct computations given in \cite{Kaw08pre}.
It is desirable to have a conceptual proof of Theorem \ref{thm:e_1}. \par

\section{Proof of Proposition \ref{prop:a3}}
\label{sec:a3}

In this section we prove Proposition \ref{prop:a3} based on the following theorem of Harris.
\begin{theorem}[\cite{Harris}, Theorem 6.5]\label{thm:Harris} For any $g \geq 3$ there exists 
a hyperelliptic Riemann surface $C_0$ such that the inclusion \eqref{eq:hypell}
is the equality
$$
Q(U(C_0)) = (\text{the $(-1)$-eigenspace of the hyperelliptic involution in $T^*_{[C_0]}\Mg$}).
$$
\end{theorem}
\noindent
In other words, the extended Johnson homomorphism $\tilde k = -[Q+\overline{Q}]$ at $[C_0]$ is non-degenerate in the normal direction to $\mathbb{H}_g$ in $\Mg$. \par
The subspace $U(C) \subset \Lambda^3(H'\oplus H'') = \bigoplus_{p+q=3}\Lambda^pH'\otimes \Lambda^qH''$ is decomposed into
$$
U(C) = \bigoplus_{p+q=3} U(C)^{p,q}, \quad
U(C)^{p,q} := U(C) \cap (\Lambda^pH'\otimes \Lambda^qH'').
$$
From the definition of the map $Q$ \eqref{eq:Q} we have 
\begin{equation}
Q\vert_{U(C)^{p,q}} = 0
\end{equation}
if $(p, q) \neq (2,1)$. Through the pairing $\langle-, -\rangle$ in \eqref{eq:M^0}, 
$U(C)^{p,q}$ and $U(C)^{3-p,3-q}$ are dual to each other. 
Hence the image of the transpose $Q: T_{[C]}\Mg \to U(C)$ is included 
in $(U(C)^{2,1})^* =  U(C)^{1,2}$, i.e., we have $
Q(T_{[C]}\Mg) \subset U(C)^{1,2}$. From the positivity of the Hermitian pairing 
\eqref{eq:pair}, we have
\begin{equation}\label{eq:pos}
\forall v \in U(C)^{1,2} \setminus\{0\}, \quad -\sqrt{-1}\langle v, \overline{v}\rangle > 0.
\end{equation}
In fact, $\{\psi_i\wedge\overline{\psi_j}\wedge\overline{\psi_k}; 1 \leq i, j, k \leq g, j< k\}$ is a orthogonal basis of $H'\otimes \Lambda^2H''(\supset U(C)^{1,2})$.\par
Now we go back to the hyperelliptic Riemann surface $[C_0] \in \mathbb{M}_3$ 
given by Theorem \ref{thm:Harris} of Harris. We choose a marking $m_0$ of $C_0$ to get a point $[C_0, m_0]$ in the Teichm\"uller space $\mathcal{T}_3$ of genus $3$.
Its Jacobian $\Jac(C_0)$ with the marking $m_0$ defines a point $[\Jac(C_0), m_0] 
\in \mathfrak{H}_3$. Then there are complex coordinates 
$(z_1, z_2, z_3)$ of $\mathcal{T}_3$ centered at $[C_0, m_0]$, and
$(w_1, w_2, w_3)$ of the Siegel upper half space $\mathfrak{H}_3$ 
centered at 
$[\Jac(C_0), m_0]$ such that the locus $\{z_1=0\}$ coincides with 
the inverse image of $\mathbb{H}_3$ locally, and that 
$$
\Jac(z_1, z_2, z_3) = ({z_1}^2, z_2, z_3) = (w_1, w_2, w_3),
$$
since the Jacobian map $\Jac$ is a $2$ to $1$ map in the normal direction of the submanifold $\mathbb{H}_3$. There are $U(C)^{1,2}$-valued functions $v_1, v_2, v_3$ defined around  
$[C_0, m_0] \in \mathcal{T}_3$ such that $\tilde k = v_1 dz_1+v_2dz_2+v_3dz_3$
on $T\mathcal{T}_3$. Then we have $v_1(0, *, *) \neq 0$ from Theorem \ref{thm:Harris}, $v_2(0, *, *) = v_3(0, *, *) = 0$ from \eqref{eq:hypell}, and
the differential form $e_1^J$ equals 
$$
\sum^3_{i,j=1}\langle v_i, \overline{v_j}\rangle 
%- \overline{\langle v_i, \overline{v_j}\rangle}
dz_i\wedge d\overline{z_j}
$$
up to non-zero constant factor. Since $\frac{\pa}{\pa w_1} = 
\frac1{2z_1}\frac{\pa}{\pa z_1}$, the amount $z_1\cdot \Lambda e_1^J$ 
is a combination of the amounts
$$
\frac1{z_1}\langle v_1, \overline{v_1}\rangle, \,\,
\langle v_1, \overline{v_j}\rangle, \,\,
\langle v_i, \overline{v_1}\rangle \,\,\text{and}\,\,
z_1\langle v_i, \overline{v_j}\rangle\quad (i, j = 2, 3)
$$
%and their complex conjugate 
with coefficients in non-zero locally bounded functions near $[C_0, m_0]$.
The amounts $\langle v_1, \overline{v_j}\rangle$, 
$\langle v_i, \overline{v_1}\rangle$ and 
$z_1\langle v_i, \overline{v_j}\rangle$ go to $0$ as $z_1 \to 0$. 
Together with \eqref{eq:pos}, this implies
$$
|z_1|\cdot |\Lambda e_1^J| \geq \frac{(\text{positive constant})}{|z_1|}
|\langle v_1, \overline{v_1}\rangle| \longrightarrow +\infty%\quad \text{as $z_1 \to 0$}.
$$
as $z_1 \to 0$. 
This completes the proof of Proposition \ref{prop:a3} and that of Theorem \ref{thm:g3}. \qed\par
\bigskip
%%%% 参考文献

\vskip 10mm
\noindent
Department of Mathematical Sciences, \\
University of Tokyo \\
3-8-1 Komaba, Meguro-ku, Tokyo, \\
153-8914, JAPAN. \\
kawazumi@ms.u-tokyo.ac.jp\\
\\
\end{document}